# $m - Magic$ Labeling on Anti Fuzzy Path Graph

Yeni Rahma Oktaviani[1, a)], Toto Nusantara[2, b)], Santi Irawati[2)]

[1]*Graduate Student of Department of Mathematics, FMIPA Universitas Negeri Malang, Indonesia,*
[2] *Department of Mathematics, FMIPA Universitas Negeri Malang, Indonesia,*

[b)] *Corresponding author: toto.nusantara.fmipa@um.ac.id*
[a)] *yeni.rahma.2103138@students.um.ac.id*

***Abstract.*** The topic of $m - magic$ labeling is discussed in this article. The application of these problem is applied to path anti fuzzy and path anti fuzzy bipolar. Anti fuzzy and anti fuzzy bipolar graphs are new concepts of fuzzy graphs. This article, to discover $m - magic$ anti fuzzy and anti fuzzy bipolar graphs through the adaptation of some previous studies on labeling magic, bi-magic, tri-magic anti fuzzy graphs, and anti fuzzy bipolar graphs. The graph path is the subject of $m - magic$ anti fuzzy and anti fuzzy bipolar graph study. The results of this study show that, for every natural number $n$ dan $n \geq (2m + 1 + ma)$ with $a = \{0,1,2,3, ..., l\}$, there is an anti fuzzy $m - magic$ graph path labeling. While in $m > 2$ in terms of 2 cases. Case 1 when for odd $n$ with $n \equiv 2m + 1 \bmod (m)$ and for even $n$ with $n \equiv 2m + 1 \bmod (2m)$, case 2 when $n \equiv m + 1 \bmod (2m)$ for $m$ odd and $n \neq m + 1$ there is $m - magic$ graph path anti fuzzy bipolar labeling.

## INTRODUCTION

Zadeh [1] first introduce the mathematical concept of uncertainty commonly known as fuzzy. Further developments regarding fuzzy graphs were later introduced by Rosenfeld [2]. Some of these fundamental ideas in graph theory, like paths and cycles, have been given fuzzy analogies. Fuzzy graphs are generalized form of crisp graphs, therefore problems that cannot be solved with crisp graphs can be solved with the help of fuzzy graphs.

The bipolar fuzzy graph set, an extension of the fuzzy graph set, was presented by Zhang [3]. Degrees of participation in the collection of bipolar fuzzy graphs are between [-1,1]. A membership degree of 0 on a member means that the member does not meet the nature of bipolar fuzzy set, while the membership degrees $[-1,0)$ and $(0,1]$ of a member indicate that the member satisfies the nature of the bipolar fuzzy set. Although the bipolar fuzzy set and the fuzzy set appear to be the same, they are actually two separate sets.

The generalization of fuzzy graphs is an anti fuzzy graph which is a new concept and was introduced by Akram [4]. Sanggoor & Alwan [5] are among those who have worked on the creation of anti fuzzy graphs with reference to the notation of anti fuzzy graphs, some anti fuzzy graph features and operations, and their application to regular and irregular graphs. Further studies related to anti fuzzy graphs were also conducted by Muthuraj & Sasireka [6] on the concept of anti fuzzy graphs and the characteristics of regular and irregular degrees of anti fuzzy graphs. Anti fuzzy graphs and fuzzy graphs are different, it is found in the relationship between the two. The anti fuzzy graph $G = (V, \sigma, \mu)$, where $\sigma: V \to [0,1]$ is the set in $V$ and $\mu: V \times V \to [0,1]$ satisfies the relation when $\mu(vu) \geq \sigma(v) \vee \sigma(u)$ for all $vu \in E$. About the fuzzy graph is relation, it is $\mu(vu) \leq \sigma(v) \vee \sigma(u)$. The bipolar anti fuzzy graph can be defined using the analogy process between the bipolar fuzzy graph and the anti fuzzy graph.

The magic labeling of fuzzy graphs can be seen in [7,8]. The labeling of a fuzzy graph $G = (V, \sigma, \mu)$ is said to be a magic fuzzy graph if $\sigma(v) + \mu(vu) + \sigma(u) = k$ for $v, u \in V$ with $k$ constant and represented by $m_0(G)$. According to Sheeba's research [7], there is a special labeling for butterfly graphs, pan graphs, wheel graphs, helm graphs, and bull graphs in fuzzy graphs. According to research by Nagorgani, et al [8], fuzzy graphs with paths, cycles, and stars have a magic labeling. The similar study was done on the magic labeling of anti fuzzy graphs on path, star, and cycle graphs by Brata, et al [14]. The magic labeling of anti fuzzy bipolar graphs on graphs path, stars, and cycles was also covered in study by Firmansa, et al [15].

The bi-magic labeling notation, which includes a value of 2 magic, was introduce by Babujee [9]. The labeling of a fuzzy graph $G = (V, \sigma, \mu)$ is said to be a bi-magic fuzzy graph if $\sigma(v) + \mu(vu) + \sigma(u) = k_1$ or $k_2$ for $v, u \in V$ with $k_1, k_2$ different constants and symbolized $Bm_0(G)$. The existence of bi-magic labeling of fuzzy bistar graphs and magic labeling of star fuzzy graphs has been demonstrated in research on the bi-magic labeling of fuzzy graphs, among others, by Thirusangu & Jeevitha [10]. Further research by Bibi & Devi [11] found bi-magic labeling of fuzzy graphs, cycles, and stars.

In Sumathi & Monigeetha's research on isomorphic 5-diameter caterpillar graphs and jahangir graphs [12,13]. The labeling of a fuzzy graph $G = (V, \sigma, \mu)$ is said to be a tri-magic fuzzy graph if $\sigma(v) + \mu(vu) + \sigma(u) = k_1$ or $k_2$ or $k_3$ for $v, u \in V$ with $k_1, k_2, k_3$ different constants and represented $Tm_0(G)$. Various cases of labeling magic, bi-magic, and fuzzy tri-magic graphs can be analogous to the labeling of magic, bi-magic, and tri-magic anti fuzzy graphs. Based on the analogy process of the various labeling cases above, the labeling of fuzzy graphs $G = (V, \sigma, \mu)$ is said to be a fuzzy $m - magic$ graph if $\sigma(v) + \mu(vu) + \sigma(u) = k_1$ or $k_2$ or $k_3, \dots, k_m$ for $v, u \in V$ with $k_1, k_2, k_3, \dots, k_m$ is a different constant and is represented $M_m(G)$.

According to Brata's research [14] on magic, bi-magic, and tri-magic labeling on anti fuzzy graphs, among other things, these is a magic labeling of anti fuzzy graph path and stars for any natural number $n$ with $m = 1$, and there is an anti fuzzy magic labeling of cycle graphs for odd $n$ with $n \geq 3$. Moreover, Brata discovered bi-magic labeling of anti fuzzy graphs on the bistar, star, and path with $m = 2$. The cases on magic labeling and bi-magic anti fuzzy graphs can be applied to magic labeling, bi-magic anti fuzzy bipolar graphs. According to Firmansa [15], anti fuzzy bipolar graphs bi-magic labeled to paths, cycles, and stars with $m = 2$. Therefore, based on the analogy process, the above case raises the question of whether there is a labeling of $m - magic$ anti fuzzy graphs dan $m - magic$ anti fuzzy bipolar graphs with $m > 2$.

In this study, $m - magic$ graph path labeling anti fuzzy and anti fuzzy bipolar with $m > 2$ will be proven. The research in this article is divided into several parts, including the initial part includes the explanation of anti fuzzy and anti fuzzy bipolar graphs and the magic labeling of anti fuzzy and anti fuzzy bipolar graphs. The results section presents the presence of $m - magic$ graph path labeling, anti fuzzy and anti fuzzy bipolar. The final section will present the conclusion.

## PRELIMINARIES

Definitions pertaining to the research setting will be provide in this part, specifically for $m - magic$ labeling anti fuzzy bipolar graphs path, as follows :

**Definition 1**. Let crisp graph $G^* = (V, E)$ with $V \neq \emptyset$ and $E \subseteq V \times V$. Anti fuzzy graph $G = (V, \sigma, \mu)$ is the pair of functions $\sigma: V \to [0,1]$ and $\mu: V \times V \to [0,1]$, where for all $v, u \in V$, satisfy $\mu(v, u) \geq \sigma(v) \vee \sigma(u)$ and $\mu$ is a fuzzy relation in $\sigma$.

According on definition 1 above, operation ($\vee$) means maximum or can be written $\sigma(v) \vee \sigma(u) = \max(\sigma(v), \sigma(u))$.

**Example 1.** Crisp graph $G^* = (V, E)$ with $V = \{v_1, v_2, v_3\}$ and $E = \{v_1v_2, v_2v_3, v_3v_1\}$. Suppose $\sigma$ the fuzzy set in $V$ and $\mu$ the anti fuzzy function in $E \subseteq V \times V$ is defined as follows

|   | $v_1$ | $v_2$ | $v_3$ |
|---|---|---|---|
| $\sigma$ | 0,2 | 0,3 | 0,5 |

|   | $v_1v_2$ | $v_2v_3$ | $v_3v_1$ |
|---|---|---|---|
| $\mu$ | 0,8 | 0,8 | 0,6 |

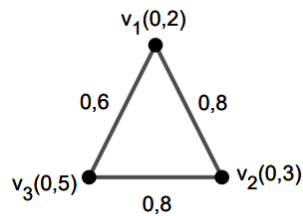

Figure 1. Anti fuzzy graph

**Definition 2**. Bipolar anti fuzzy graph with set $V$ defined as pairs $G = (A, B)$ where $A = (\sigma^P, \sigma^N)$ the set of anti fuzzy bipolar in $V$ and $B = (\mu^P, \mu^N)$ the set of anti fuzzy bipolar in $E \subseteq V \times V$, for all $vu \in E$ satisfies $\mu^P(vu) \geq \sigma^P(v) \vee \sigma^P(u)$ and $\mu^N(vu) \leq \sigma^N(v) \wedge \sigma^N(u)$

Definition 2 is obtained through the analogy of concepts from the anti fuzzy graph and the bipolar fuzzy set [3]. Operation $\wedge$ means minimum.

**Example 2.** Crisp graph $G^* = (V, E)$ with $V = \{v, u_1, u_2, u_3\}$ and $E = \{vu_1, vu_2, vu_3\}$. Let $A = (\sigma^P, \sigma^N)$ the set of anti fuzzy bipolar in $V$ and $B = (\mu^P, \mu^N)$ the set of anti fuzzy bipolar in $E \subseteq V \times V$ are defined as follows

|            | $v$  | $u_1$ | $u_2$ | $u_3$ |
|------------|------|-------|-------|-------|
| $\sigma^P$ | 0,1  | 0,2   | 0,3   | 0,4   |
| $\sigma^N$ | -0,5 | -0,6  | -0,7  | -0,8  |

|          | $vu_1$ | $vu_2$ | $vu_3$ |
|----------|--------|--------|--------|
| $\mu^P$  | 0,3    | 0,4    | 0,5    |
| $\mu^N$  | -0,7   | -0,8   | -0,9   |

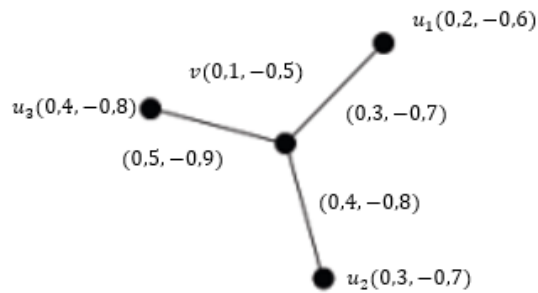

Figure 2. Bipolar anti fuzzy graph

**Definition 3.** [16] Graph path is a graph with a track between the point and the side traversed by which nothing repeats. Graph path have $n$ vertice and $n - 1$ edge are denoted with $P_n$.

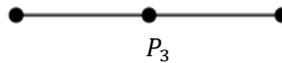

Figure 3. Graph path $P_3$

Furthermore, anti fuzzy path graph and anti fuzzy bipolar graph path are defined as follows

**Definition 4.** Graph path $P_n$ be said to be anti fuzzy graph path if $\mu(v_i v_{i+1}) \geq \sigma(v_i) \vee \sigma(v_{i+1})$
**Definition 5.** Graph path $P_n$ is said to be anti fuzzy bipolar graph path if it is $\mu^P(v_i v_{i+1}) \geq \sigma^P(v_i) \vee \sigma^P(v_{i+1})$ and $\mu^N(v_i v_{i+1}) \geq \sigma^N(v_i) \vee \sigma^N(v_{i+1})$ for all $v_i v_{i+1} \in E$
**Definition 6.** [7,8] Labeling a fuzzy graph $G = (V, \sigma, \mu)$ is said to be a magic fuzzy graph if $\sigma(v) + \mu(vu) + \sigma(u) = k$ for $v, u \in V$ with $k$ constant and represented by $m_0(G)$

Brata [14] performed the following magic labeling of anti fuzzy graphs path in Theorem 1

**Theorem 1.** [14] Let $P_n$, $n > 1$ anti fuzzy graph path. Then there is the magic labeling $P_n$

If $n = 1$ in Theorem 1 does not form a graph path due to the ansence of edges, $n = 2$ does not yet produce the value of magic constant where there is just one graph edge. Hence, the magic labeling of anti fuzzy graphs path is satisfied by $n \geq 3$. Theorem 1 results in the magic constant $m_0(G) = 3nd$, with the vertices labeled as $\sigma(v_i) = id$ and $\sigma(v_{i+1}) = (i + 1)$, and the edges labeled as $\mu(v_i v_{i+1})$ which is $(3n - 2i - 1)d$. The vertices should be labeled in order while being labeled.

**Example 3.** Crisp graph $G^* = (V, E)$ with $V = \{v_1, v_2, v_3, v_4, v_5\}$ and $E = \{v_1 v_2, v_2 v_3, v_3 v_4, v_4 v_5\}$. Let $\sigma$ the fuzzy set in $V$ and the $\mu$ of the anti fuzzy function in $E \subseteq V \times V$ are defined as follows

|          | $v_1$ | $v_2$ | $v_3$ | $v_4$ | $v_5$ |
|----------|-------|-------|-------|-------|-------|
| $\sigma$ | 0,01  | 0,02  | 0,03  | 0,04  | 0,05  |

|       | $v_1 v_2$ | $v_2 v_3$ | $v_3 v_4$ | $v_4 v_5$ |
|-------|-----------|-----------|-----------|-----------|
| $\mu$ | 0,12      | 0,1       | 0,08      | 0,06      |

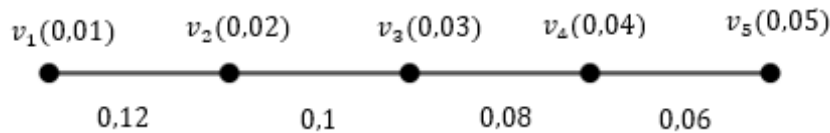

Figure 4. Magic labeling of anti fuzzy path graph $P_4$

In Figure 4 the labeling has a magic constant of $k_1 = 0,15$, so the labeling is a magic labeling of anti fuzzy graph path.

**Definition 7**. [9] The labeling of a fuzzy graph $G = (V, \sigma, \mu)$ is said to be a bi-magic fuzzy graph if $\sigma(v) + \mu(vu) + \sigma(u) = k_1$ or $k_2$ for $v, u \in V$ with $k_1, k_2$ different constants and symbolized $Bm_0(G)$

Brata [14] performed the following bi-magic labeling of anti fuzzy graphs path in Theorem 2.

**Theorem 2.** For all a $n \geq 5$ and odd $n$, the graph path $P_n$ satisfies the bi-magic labeling anti fuzzy graph

Because there are no edges, $n = 1$ does not form a graph path in Theorem 2, and $n = 2$ does not yet create the bi-magic constant value for graphs with a single edge. In the case of $n = 3$, there are two edges that are unable to determine the bi-magic constant's value, and in the case of $n = 4$, there are theree edges that do not split the number of bi-magic constants equally. As a result, the bi-magic labeling of anti fuzzy graphs path is satisfied for $n = 5$ and odd $n$. Theorem 2 yields bi-magic constants $Bm_1(G) = (2n + 2)d$ and $Bm_2(G) = (2n + 7)d$ where the edge labeling is $\mu(v_i v_{i+1}) = (2n + 1 - 2i)d$ for $Bm_1(G)$, $\mu(v_i v_{i+1}) = (2n + 6 - 2i)d$ for $Bm_2(G)$, the vertice labeling is $\sigma(v_i) = id$ and $\sigma(v_{i+1}) = (i + 1)$. The vertices should be labeled in order while being labeled. Because Theorem 2's labeling does not satisfies the labeling of an edge in an anti fuzzy network, which is $\mu(v_i v_{i+1}) \geq \sigma(v_i) \vee \sigma(v_{i+1})$, it cannot be applied to tri-magic labeling.

**Definition 8**. [12,13] The labeling of fuzzy graph $G = (V, \sigma, \mu)$ is said to be a tri-magic fuzzy graph if $\sigma(v) + \mu(vu) + \sigma(u) = k_1$ or $k_2$ or $k_3$ for $v, u \in V$ with $k_1, k_2, k_3$ different constants and symbolized $Tm_0(G)$

The authors started looking at the tri-magic instance since they were unable to locate the graph of the anti fuzzy path on the labeling of tri-magic. Because labeling on the edge will not satisfy the idea of anti fuzzy graph labeling, magic and bi-magic edge labeling cannot be applied to tri-magic, although the process for labeling them is the same.

The idea of labeling $m - magic$ fuzzy graphs is taken from those definitions of labeling $m - magic$ fuzzy graphs

**Definition 9**. The labeling of a fuzzy graph $G = (V, \sigma, \mu)$ is said to be $m - magic$ fuzzy graph if $\sigma(v) + \mu(vu) + \sigma(u) = k_1$ or $k_2$ or $k_3, \dots, k_m$ for $v, u \in V$ with $k_1, k_2, k_3, \dots, k_m$ is a different constant and is symbolized $M_m(G)$.

Moreover, Theorem 3 pertaining to Firmansa's [15] anti fuzzy bipolar graph's path magic labeling states the following:

**Theorem 3.** Let $P_n, n > 1$ graph path $P_n$ anti fuzzy bipolar. So there is a magic labeling on $P_n$

The bipolar anti fuzzy path's magic labeling graph, which had a length of $n > 1$, was separated into two cases for $n$ odd and $n$ even in Theorem 3. In labeling, each case has a positive and negative label. Edge labeling $\mu^P(v_i v_{i+1}) = (6n - 4i - 1)d$ for $m_0^P(G)$, and $\mu^N(v_i v_{i+1}) = (1 + 4i - 6n)d$ for $m_0^N(G)$. To the $m - magic$ labeling of bipolar anti fuzzy graph path, the magic labeling of bipolar anti fuzzy graph path can be applied. Then $m - magic$ labeling of bipolar fuzzy graph is defined as follows.

**Definition 10**. The labeling of a fuzzy graph $G = (V, \sigma, \mu)$ is said to be a bipolar fuzzy graph $m - magic$ if
i. $\sigma^P(v) + \mu^P(vu) + \sigma^P(u) = k_1^P$ or $k_2^P$ or $k_3^P, \dots, k_m^P$ for $v, u \in V$ with $k_1, k_2, k_3, \dots, k_m$ different constants and symbolized $M_m^P(G)$.
ii. $\sigma^N(v) + \mu^N(vu) + \sigma^N(u) = k_1^N$ or $k_2^N$ or $k_3^N, \dots, k_m^N$ for $v, u \in V$ with $k_1, k_2, k_3, \dots, k_m$ different constants and symbolized $M_m^N(G)$.

## RESULTS AND DISCUSSION

Theorem 4 for $m - magic$ labeling of anti fuzzy graph path, where this labeling can also apply to magic and bi-magic labeling, is where we provide the findings of our research on the magic labeling of anti fuzzy graphs in this part. Labeling $m - magic$ bipolar anti fuzzy graph path, Theorem 5. The Theorem is a novel idea that builds on earlier work on the magic labeling of anti fuzzy and anti fuzzy bipolar graphs by Brata [14] and Firmansa [15].

**Theorem 4.** Let $P_n$, $n \geq (2m + 1 + ma)$ with $a = \{0,1,2, \dots, l\}$ for $m > 2$ anti fuzzy graph path. Then there is a $m - magic$ labeling on $P_n$.

**Proof:** Let $G = (V, \sigma, \mu)$ is an anti fuzzy path graph $P_n$ with length $n \geq (2m + 1 + ma)$. Let $V = \{v_1, v_2, v_3, \ldots, v_n\}$. Define labeling as follow
1. For vertices, defines $\sigma(v_i) = id$, where $1 \leq i \leq n$
2. For edges, defines $\mu(v_i v_{i+1}) = (3n - 2i + 1)d$ where $1 \leq i \leq \frac{n-1}{m}$
3. For edges, defines $\mu(v_i v_{i+1}) = (3n - 2i + 4)d$ where $\frac{n+m-1}{m} \leq i \leq \frac{2n-2}{m}$
4. For edges, defines $\mu(v_i v_{i+1}) = (3n - 2i + 8)d$ where $\frac{2n+m-2}{m} \leq i \leq \frac{3n-3}{m}$
5. For edges, defines $\mu(v_i v_{i+1}) = (3n - 2i + 10)d$ where $\frac{3n+m-3}{m} \leq i \leq \frac{4n-4}{m}$
6. For edges, defines $\mu(v_i v_{i+1}) = (3n - 2i + 12)d$ where $\frac{4n+m-4}{m} \leq i \leq \frac{5n-5}{m}$
   $\vdots$
7. For edges, defines $\mu(v_i v_{i+1}) = (3n - 2i + (2m + 2))d$ where $\frac{kn+m-k}{m} \leq i \leq \frac{kn-k}{m}$ and $d$ is given by

$$d = \begin{cases} 10^{-2}, (2m+1) \leq n < 31 \\ 10^{-3}, 31 \leq n < 331 \\ 10^{-(j+4)}, 331 \times 10^j \leq n \leq 331 \times 10^{j+1} \end{cases}$$

Notes: $j = million$

Thus we obtained the $m - magic$ constant is, $1 \leq i \leq n-1$

$M_1(G) = \sigma(v_i) + \mu(v_i v_{i+1}) + \sigma(v_{i+1})$
$= (i)d + (3n + 1 - 2i)d + (i + 1)d$
$= (3n + 2)d$

$M_2(G) = \sigma(v_i) + \mu(v_i v_{i+1}) + \sigma(v_{i+1})$
$= (i)d + (3n + 4 - 2i)d + (i + 1)d$
$= (3n + 5)d$

$M_3(G) = \sigma(v_i) + \mu(v_i v_{i+1}) + \sigma(v_{i+1})$
$= (i)d + (3n + 8 - 2i)d + (i + 1)d$
$= (3n + 9)d$

$M_4(G) = \sigma(v_i) + \mu(v_i v_{i+1}) + \sigma(v_{i+1})$
$= (i)d + (3n + 10 - 2i)d + (i + 1)d$
$= (3n + 11)d$

$M_5(G) = \sigma(v_i) + \mu(v_i v_{i+1}) + \sigma(v_{i+1})$
$= (i)d + (3n + 12 - 2i)d + (i + 1)d$
$= (3n + 13)d$
$\vdots$

$M_m(G) = \sigma(v_i) + \mu(v_i v_{i+1}) + \sigma(v_{i+1})$
$= (i)d + (3n + 2m + 2 - 2i)d + (i + 1)d$
$= (3n + 2m + 3)d$

It is sufficient to prove $\mu(v_i v_{i+1}) \geq \max(\sigma(v_i), \sigma(v_{i+1})), 1 \leq i \leq n - 1$
$\mu(v_i v_{i+1}) = (3n - 2i + (2m + 2))d$
$\geq \max((i)d, (i+1)d)$
$\geq \max(\sigma(v_i), \sigma(v_{i+1}))$

The following is given an example of $m - magic$ labeling of anti fuzzy graph path

**Example 4.** Crisp graph $G^* = (V, E)$ with $V = \{v_1, v_2, v_3, \ldots, v_9\}$ and $E = \{v_1 v_2, v_2 v_3, v_3 v_4, \ldots, v_8 v_9\}$. Let $\sigma$ the fuzzy set in $V$ and $\mu$ of the anti fuzzy function in $E \subseteq V \times V$ are defined as follows

|   | $v_1$ | $v_2$ | $v_3$ | $v_4$ | $v_5$ | $v_6$ | $v_7$ | $v_8$ | $v_9$ |
|---|---|---|---|---|---|---|---|---|---|
| $\sigma$ | 0,01 | 0,02 | 0,03 | 0,04 | 0,05 | 0,06 | 0,07 | 0,08 | 0,09 |

|   | $v_1 v_2$ | $v_2 v_3$ | $v_3 v_4$ | $v_4 v_5$ | $v_5 v_6$ | $v_6 v_7$ | $v_7 v_8$ | $v_8 v_9$ |
|---|---|---|---|---|---|---|---|---|
| $\mu$ | 0,26 | 0,24 | 0,25 | 0,23 | 0,25 | 0,23 | 0,23 | 0,21 |

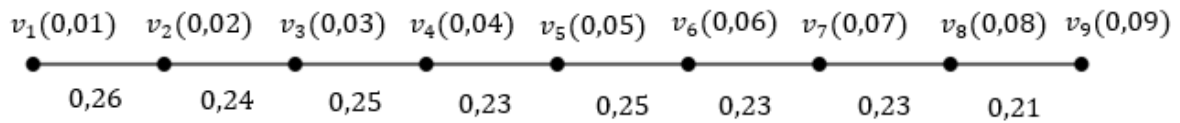

Figure 5. Anti fuzzy graph path $P_9$

In Figure 5, the labeling is an $m - magic$ labeling of an anti fuzzy graph path with $m = 4$. The labeling has four magic constants, $k_1 = 0,29$, $k_2 = 0,32$, $k_3 = 0,36$ and $k_4 = 0,36$.

**Theorem 5.** Let $P_n$, $n \geq (2m + 1 + ma)$ with $a = \{0,1,2, ..., l\}$ for $m > 2$ anti fuzzy bipolar graph path. Then there is a $m - magic$ labeling on $P_{npath}$.

**Proof:** Let $G = (V, \sigma, \mu)$ is an anti fuzzy graph bipolar $P_n$ with length $n \geq (2m + 1 + ma)$. Let $V = \{v_1, v_2, v_3, ..., v_n\}$. We see the following cases.

**Case 1:** When $n \equiv 2m + 1 \mod (m)$ for $m$ even, and when $n \equiv 2m + 1 \mod (2m)$ for $m$ odd

Define labeling as follow :

1. For odd vertices, defines $\sigma^P(v_i) = 2i - 1d$ and $\sigma^N(v_i) = (1 - 2i)d$ where $1 \leq i \leq \frac{n+1}{2}$
2. For even vertices, defines $\sigma^P(v_i) = (2i)d$ and $\sigma^N(v_i) = -(2i)d$ where $1 \leq i \leq \frac{n-1}{2}$
3. For edges, defines $\mu^P(v_i v_{i+1}) = (6n - 4i - 1)d$ and $\mu^N(v_i v_{i+1}) = (1 + 4i - 6n)d$ where $1 \leq i \leq \frac{n-1}{m}$
4. For edges, defines $\mu^P(v_i v_{i+1}) = (7n - 4i - 2)d$ and $\mu^N(v_i v_{i+1}) = (2 + 4i - 7n)d$ where $\frac{n+m-1}{m} \leq i \leq \frac{2n-2}{m}$
5. For edges, defines $\mu^P(v_i v_{i+1}) = (8n - 4i - 3)d$ and $\mu^N(v_i v_{i+1}) = (3 + 4i - 8n)d$ where $\frac{2n+m-2}{m} \leq i \leq \frac{3n-3}{m}$
6. For edges, defines $\mu^P(v_i v_{i+1}) = (9n - 4i - 4)d$ and $\mu^N(v_i v_{i+1}) = (4 + 4i - 9n)d$ where $\frac{3n+m-3}{m} \leq i \leq \frac{4n-4}{m}$
7. For edges, defines $\mu^P(v_i v_{i+1}) = (10n - 4i - 5)d$ and $\mu^N(v_i v_{i+1}) = (5 + 4i - 10n)d$ where $\frac{4n+m-4}{m} \leq i \leq \frac{5n-5}{m}$

   $\vdots$

8. For edges, defines $\mu^P(v_i v_{i+1}) = ((m+5)n - 4i - m)d$ and $\mu^N(v_i v_{i+1}) = (m + 4i - (m+5)n)d$ where $\frac{kn+m-k}{m} \leq i \leq \frac{kn-k}{m}$ and $d$ is given by

$$d = \begin{cases} 10^{-2}, (2m+1) \leq n < 11 \\ 10^{-3}, 11 \leq n < 35 \\ 10^4, 35 \leq 334 \\ 10^{-(j+4)}, 334 \times 10^j \leq n \leq 334 \times 10^{j+1} \end{cases}$$

Notes: $j = million$

Thus we obtained the $m - magic$ constant is , $1 \leq i \leq n - 1$

$$M_1^P(G) = \sigma^P(v_i) + \mu^P(v_i v_{i+1}) + \sigma^P(v_{i+1})$$
$$= (2i - 1)d + (6n - 4i - 1)d + (2i)d$$
$$= 6nd$$
$$M_1^N(G) = \sigma^N(v_i) + \mu^N(v_i v_{i+1}) + \sigma^N(v_{i+1})$$
$$= (1 - 2i)d + (1 + 4i - 6n)d - (2i)d$$
$$= -6nd$$
$$M_2^P(G) = \sigma^P(v_i) + \mu^P(v_i v_{i+1}) + \sigma^P(v_{i+1})$$
$$= (2i - 1)d + (7n - 4i - 2)d + (2i)d$$
$$= (7n - 1)d$$
$$M_2^N(G) = \sigma^N(v_i) + \mu^N(v_i v_{i+1}) + \sigma^N(v_{i+1})$$
$$= (1 - 2i)d + (2 + 4i - 7n)d - (2i)d$$
$$= -(7n - 1)d$$
$$M_3^P(G) = \sigma^P(v_i) + \mu^P(v_i v_{i+1}) + \sigma^P(v_{i+1})$$
$$= (2i - 1)d + (8n - 4i - 3)d + (2i)d$$
$$= (8n - 2)d$$
$$M_3^N(G) = \sigma^N(v_i) + \mu^N(v_i v_{i+1}) + \sigma^N(v_{i+1})$$
$$= (1 - 2i)d + (3 + 4i - 8n)d - (2i)d$$
$$= -(8n - 2)d$$
$$M_4^P(G) = \sigma^P(v_i) + \mu^P(v_i v_{i+1}) + \sigma^P(v_{i+1})$$

$$\begin{aligned}
&= (2i-1)d + (9n-4i-4)d + (2i)d \\
&= (9n-3)d \\
M_4^N(G) &= \sigma^N(v_i) + \mu^N(v_i v_{i+1}) + \sigma^N(v_{i+1}) \\
&= (2i-1)d + (4-4i-9n)d - (2i)d \\
&= -(9n-3)d \\
M_5^P(G) &= \sigma^P(v_i) + \mu^P(v_i v_{i+1}) + \sigma^P(v_{i+1}) \\
&= (2i-1)d + (10n-4i-5)d + (2i)d \\
&= (10n-4)d \\
M_5^N(G) &= \sigma^N(v_i) + \mu^N(v_i v_{i+1}) + \sigma^N(v_{i+1}) \\
&= (1-2i)d + (5+4i-10n)d - (2i)d \\
&= -(10n-4)d \\
&\vdots \\
M_m^P(G) &= \sigma^P(v_i) + \mu^P(v_i v_{i+1}) + \sigma^P(v_{i+1}) \\
&= (2i-1)d + ((m+5)n - 4i - m)d + (2i)d \\
&= ((m+5)n - (m-1))d \\
M_m^N(G) &= \sigma^N(v_i) + \mu^N(v_i v_{i+1}) + \sigma^N(v_{i+1}) \\
&= (1-2i)d + (m+4i-(m+5)n)d - (2i)d \\
&= -((m+5)n - (m-1))d
\end{aligned}$$

It is sufficient to prove $\mu^P(v_i v_{i+1}) \geq \max(\sigma^P(v_i), \sigma^P(v_{i+1})), 1 \leq i \leq n-1$

$$\begin{aligned}
\mu^P(v_i v_{i+1}) &= ((m+5)n - (m-1))d \\
&\geq \max((2i-1)d, (2i)d) \\
&\geq \max(\sigma^P(v_i), \sigma^P(v_{i+1}))
\end{aligned}$$

And $\mu^N(v_i v_{i+1}) \leq \max(\sigma^N(v_i), \sigma^N(v_{i+1})), 1 \leq i \leq n-1$

$$\begin{aligned}
\mu^N(v_i v_{i+1}) &= -((m+5)n - (m-1))d \\
&\leq \max((1-2i)d, (-2i)d) \\
&\leq \max(\sigma^N(v_i), \sigma^N(v_{i+1}))
\end{aligned}$$

**Case 2:** When $n \equiv m+1 \, mod \, (2m)$ for $m$ odd and $n \neq m+1$

Define labeling as follow

1. For odd vertices, defines $\sigma^p(v_i) = 2i - 1d$ and $\sigma^N(v_i) = (1-2i)d$ where $1 \leq i \leq \frac{n}{2}$
2. For even vertices, defines $\sigma^p(v_i) = (2i)d$ and $\sigma^N(v_i) = -(2i)d$ where $1 \leq i \leq \frac{n}{2}$
3. For edges, defines $\mu^P(v_i v_{i+1}) = (6n - 4i - 1)d$ and $\mu^N(v_i v_{i+1}) = (1 + 4i - 6n)d$ where $1 \leq i \leq \frac{n-1}{m}$
4. For edges, defines $\mu^P(v_i v_{i+1}) = (7n - 4i - 2)d$ and $\mu^N(v_i v_{i+1}) = (2 + 4i - 7n)d$ where $\frac{n+m-1}{m} \leq i \leq \frac{2n-2}{m}$
5. For edges, defines $\mu^P(v_i v_{i+1}) = (8n - 4i - 3)d$ and $\mu^N(v_i v_{i+1}) = (3 + 4i - 8n)d$ where $\frac{2n+m-2}{m} \leq i \leq \frac{3n-3}{m}$
6. For edges, defines $\mu^P(v_i v_{i+1}) = (9n - 4i - 4)d$ and $\mu^N(v_i v_{i+1}) = (4 + 4i - 9n)d$ where $\frac{3n+m-3}{m} \leq i \leq \frac{4n-4}{m}$
7. For edges, defines $\mu^P(v_i v_{i+1}) = (10n - 4i - 5)d$ and $\mu^N(v_i v_{i+1}) = (5 + 4i - 10n)d$ where $\frac{4n+m-4}{m} \leq i \leq \frac{5n-5}{m}$

$$\vdots \qquad \vdots \qquad \vdots$$

8. For edges, defines $\mu^P(v_i v_{i+1}) = ((m+5)n - 4i - m)d$ and $\mu^N(v_i v_{i+1}) = (m + 4i - (m+5)n)d$ where $\frac{kn+m-k}{m} \leq i \leq \frac{kn-k}{m}$ and $d$ is given by

$$d = \begin{cases} 10^{-2}, (2m+1) \leq n < 11 \\ 10^{-3}, 11 \leq n < 35 \\ 10^4, 35 \leq n \leq 334 \\ 10^{-(j+4)}, 334 \times 10^j \leq n \leq 334 \times 10^{j+1} \end{cases}$$

Notes: $j = million$

Thus we obtained the $m - magic$ constant is , $1 \leq i \leq n-1$

$$M_1^P(G) = \sigma^P(v_i) + \mu^P(v_i v_{i+1}) + \sigma^P(v_{i+1})$$
$$= (2i-1)d + (6n-4i-1)d + (2i)d$$
$$= 6nd$$
$$M_1^N(G) = \sigma^N(v_i) + \mu^N(v_i v_{i+1}) + \sigma^N(v_{i+1})$$
$$= (1-2i)d + (1+4i-6n)d - (2i)d$$
$$= -6nd$$
$$M_2^P(G) = \sigma^P(v_i) + \mu^P(v_i v_{i+1}) + \sigma^P(v_{i+1})$$
$$= (2i-1)d + (7n-4i-2)d + (2i)d$$
$$= (7n-1)d$$
$$M_2^N(G) = \sigma^N(v_i) + \mu^N(v_i v_{i+1}) + \sigma^N(v_{i+1})$$
$$= (1-2i)d + (2+4i-7n)d - (2i)d$$
$$= -(7n-1)d$$
$$M_3^P(G) = \sigma^P(v_i) + \mu^P(v_i v_{i+1}) + \sigma^P(v_{i+1})$$
$$= (2i-1)d + (8n-4i-3)d + (2i)d$$
$$= (8n-2)d$$
$$M_3^N(G) = \sigma^N(v_i) + \mu^N(v_i v_{i+1}) + \sigma^N(v_{i+1})$$
$$= (1-2i)d + (3+4i-8n)d - (2i)d$$
$$= -(8n-2)d$$
$$M_4^P(G) = \sigma^P(v_i) + \mu^P(v_i v_{i+1}) + \sigma^P(v_{i+1})$$
$$= (2i-1)d + (9n-4i-4)d + (2i)d$$
$$= (9n-3)d$$
$$M_4^N(G) = \sigma^N(v_i) + \mu^N(v_i v_{i+1}) + \sigma^N(v_{i+1})$$
$$= (2i-1)d + (4-4i-9n)d - (2i)d$$
$$= -(9n-3)d$$
$$M_5^P(G) = \sigma^P(v_i) + \mu^P(v_i v_{i+1}) + \sigma^P(v_{i+1})$$
$$= (2i-1)d + (10n-4i-5)d + (2i)d$$
$$= (10n-4)d$$
$$M_5^N(G) = \sigma^N(v_i) + \mu^N(v_i v_{i+1}) + \sigma^N(v_{i+1})$$
$$= (1-2i)d + (5+4i-10n)d - (2i)d$$
$$= -(10n-4)d$$
$$\vdots$$
$$M_m^P(G) = \sigma^P(v_i) + \mu^P(v_i v_{i+1}) + \sigma^P(v_{i+1})$$
$$= (2i-1)d + ((m+5)n-4i-m)d + (2i)d$$
$$= ((m+5)n-(m-1))d$$
$$M_m^N(G) = \sigma^N(v_i) + \mu^N(v_i v_{i+1}) + \sigma^N(v_{i+1})$$
$$= (1-2i)d + (m+4i-(m+5)n)d - (2i)d$$
$$= -((m+5)n-(m-1))d$$

It is sufficient to prove $\mu^P(v_i v_{i+1}) \geq \max(\sigma^P(v_i), \sigma^P(v_{i+1})), 1 \leq i \leq n-1$
$$\mu^P(v_i v_{i+1}) = ((m+5)n-(m-1))d$$
$$\geq \max((2i-1)d, (2i)d)$$
$$\geq \max(\sigma^P(v_i), \sigma^P(v_{i+1}))$$
And $\mu^N(v_i v_{i+1}) \leq \max(\sigma^N(v_i), \sigma^N(v_{i+1})), 1 \leq i \leq n-1$
$$\mu^N(v_i v_{i+1}) = -((m+5)n-(m-1))d$$
$$\leq \max((1-2i)d, (-2i)d)$$
$$\leq \max(\sigma^N(v_i), \sigma^N(v_{i+1}))$$

The following is given an example of $m - magic$ labeling of anti fuzzy bipolar graph path

**Example 5.** Crisp graph $G^* = (V, E)$ with $V = \{v_1, v_2, v_3, v_4, \ldots, v_9\}$ and $E = \{v_1 v_2, v_2 v_3, v_3 v_4, \ldots, v_8 v_9\}$. Let $A = (\sigma^P, \sigma^N)$ the anti fuzzy bipolar set in $V$ and $B = (\mu^P, \mu^N)$ of the anti fuzzy bipolar function in $E \subseteq V \times V$ are defined as follows

|            | $v_1$  | $v_2$  | $v_3$  | $v_4$  | $v_5$  | $v_6$  | $v_7$  | $v_8$  | $v_9$  |
|------------|--------|--------|--------|--------|--------|--------|--------|--------|--------|
| $\sigma^P$ | 0,01   | 0,04   | 0,05   | 0,08   | 0,09   | 0,12   | 0,13   | 0,16   | 0,17   |
| $\sigma^N$ | -0,01  | -0,04  | -0,05  | -0,08  | -0,09  | -0,12  | -0,13  | -0,16  | -0,17  |

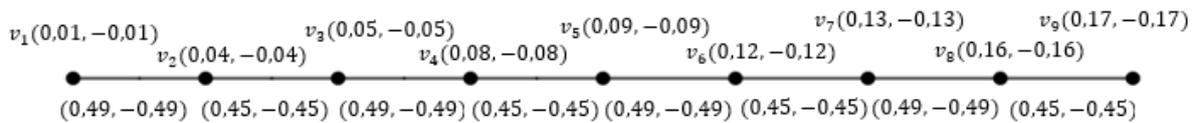

Figure 6. Anti fuzzy bipolar graph $P_9$

In Figure 6 the labeling has four magic constants, namely $k_1 = (0,54, -0,54)$ $k_2 = (0,62, -0,62)$, $k_3 = (0,7, -0,7)$ and $k_4 = (0,78, -0,78)$ so that the labeling is an $m - magic$ labeling of an anti fuzzy bipolar graph path with $m = 4$.

## CONCLUSION

In this article, we cover the labeling of $m - magic$ anti fuzzy graphs and anti fuzzy bipolar graphs. Let $G = (V, \sigma, \mu)$ is an anti fuzzy graph path $P_n$ with length of $n \geq (2m + 1 + ma)$, then there is a $m - magic$ labeling $M_m(G) = (3n + 2m + 3)d$. furthermore, let $G = (V, \sigma, \mu)$ is a graph path of anti fuzzy bipolar $P_n$ with length $n \geq (2m + 1 + ma)$, then there is a $m - magic$ labeling $M_m^P(G) = ((m + 5)n - (m - 1))d$ and $M_m^N(G) = -((m + 5)n - (m - 1))d$. Research related to magic labeling, bi-magic anti fuzzy graphs, anti fuzzy bipolar for simple graphs has been found [14,15]. Therefore, the labeling of $m - magic$ anti fuzzy graph cycle, star, bistar and labeling $m - ajaib$ anti fuzzy bipolar graph of cycle, star, bistar are still open to research.

## ACKNOWLEDGMENTS

We appreciate Ainin Yusri Saputri's assistance with this study. This article has been accepted in Proceedings of IcoMathApp 2022.